\documentclass[10pt]{article}
\usepackage{latexsym,amssymb,amsmath}

\def\N{\mathbb{N}}

\def\Z{\mathbb{Z}}
\def\R{\mathbb{R}}
\def\C{\mathbb{C}}

\def\proof{\par\medskip\noindent{\em Proof. }}
\def\eproof{\hfill{$\Box$}\bigskip}

\def\tec{\hspace{-1.6mm}{\bf. }}
\def\ds{\dots}
\def\sus{\subset}
\def\al{\alpha}
\def\be{\beta}

\def\de{\delta}

\def\ep{\varepsilon}
\def\cc{\colon}

\newtheorem{thm}{Theorem}[section]
\newtheorem{prop}[thm]{Proposition}
\newtheorem{cor}[thm]{Corollary}

\newtheorem{defi}[thm]{Definition}

\begin{document}
\title{The Newton integral and the Stirling formula}
\author{Martin Klazar}

\maketitle

\begin{abstract}
We present details of logically simplest integral  sufficient for deducing the Stirling asymptotic formula for $n!$. It is the Newton integral, 
defined as the difference of values of any primitive at the endpoints of the integration interval. We review in its framework in detail two derivations of the Stirling formula. The first approximates $\sum_{i=1}^n\log i$ with an integral and the second uses the classical gamma function and a Fubini-type result. We mention  two more integral representations of $n!$.
\end{abstract}

\section{Introduction}

Asymptotic analysis, the theory and practice of asymptotic estimates for various\,---\,often discrete\,---\,quantities, belongs to the main applications of the integral 
calculus. An archetypal example is the Stirling formula $n!\sim\sqrt{2\pi n}(\frac{n}{e})^n$ where $n$ is in $\N=\{1,2,\ds\}$, $n\to\infty$, and $n!$, the factorial of $n$, is the product $1\cdot2\cdot\ldots\cdot n$ of the first $n$ positive integers and also equals the number of $n$-tuples $(a_1,a_2,\ds,a_n)$ in $\{1,2,\ds,n\}^n$ 
such that the cardinality $|\{a_1,a_2,\ds,a_n\}|=n$. In Section 3 we present two proofs of the Stirling formula by integrals. But what kind of integrals does asymptotic analysis use, or should use?

The integrals most often used are the Riemann integral $(R)\int$, the Riemann--Stieltjes integral $(RS)\int$, the Lebesgue integral $(L)\int$, 
the Cauchy integral $(C)\int$ in $\C$ (see, for example, G.\,P. Egorychev \cite{egor} and M.\,R.~Riedel \cite{ried}), and their multivariate 
versions, especially the multivariate Cauchy integral $(MC)\int=(MC)\int\int\ds\int$ 
(see, for example, B.\,D.~McKay \cite{mcka} and R. Pemantle and M.\,C. Wilson \cite{pema_wils}). We give  four expressions of $n!$ by an integral. 
The first two are 
\begin{eqnarray*}
\log(n!)&=&c+O(1/n)+(R)\int_{1/2}^{n+1/2}\log x,\mbox{ for a $c\in\R$ and all $n\in\N$, and}\\
n!&=&(R)\int_0^{+\infty}x^ne^{-x}\\
&:=&\lim_{y\to+\infty}(R)\int_0^y x^ne^{-x},\;\mbox{ for all $n\in\N_0=\{0,\,1,\,2,\,\ds\}$}\;.
\end{eqnarray*}
We obtain them below in Propositions~\ref{appr_sum_logs} and \ref{gamma}, respectively, as Newton integrals $(N)\int$. The third and fourth expression 
of $n!$ by an $\int$ are ($n\in\N$)
\begin{eqnarray*}
\frac{1}{n!}&=&\frac{1}{2\pi i}\cdot(C)\int\frac{e^z}{z^{n+1}}\;\mbox{ and}\\
n!&=&\frac{1}{(2\pi i)^n}\cdot(MC)\int\int\ds\int\frac{(z_1+z_2+\ds+z_n)^n}{(z_1z_2\ds z_n)^2}
\end{eqnarray*}
where we integrate along counter-clockwise oriented circles in $\C$, centered at the origins. We will not consider in detail these two expressions, which are 
easy to establish by the Cauchy residue theorem. Three features set apart the last formula. The integrand is a rational and not a transcendental function. It computes $n!$ combinatorially (as the number of permutations of an $n$-element set) and not arithmetically
(as the product of the first $n$ natural numbers). Finally, the first three integral expressions for $n!$ are well known, but we have not encountered the fourth one in the literature.
We wonder if there are more simple integral representations of $n!$. One can give many more not so simple relations involving integrals and factorials. For example, F. Qi and B.-N. Guo \cite{qi_guo} present many 
integral representations of the Catalan numbers $C_n=\frac{(2n)!}{(n+1)!n!}$, yielding relations like (\cite[Theorem 3]{qi_guo})
$$
(2n)!=(n+1)!\cdot n!\cdot\frac{1}{\pi}(R)\int_0^2 x^{2n}\sqrt{4-x^2}\,dx\;.
$$

Texts on asymptotic analysis, like the books N.\,G. de Bruijn \cite{debr} or P.~Flajolet and R. Sedgewick \cite{flaj_sedg}, usually do not devote much attention to the exact definition
and properties of integrals they use and take their theory for granted, which is understandable, but some books do. For example, the monograph \cite{mont_vaug} 
by H.\,L. Montgomery and R.\,C. Vaughan has an appendix on the $(RS)\int$ in which its definition and basic properties are given. In our article we want to present derivations
of the Stirling formula with all their integral details and we aim at logical simplicity. 
Thus we need a theoretically simple integral. For example, not to take the $(C)\int$ for granted and instead to develop this powerful and versatile 
integral from scratch is not a straightforward task. It is not enough to open some of many textbooks on complex analysis because they all reach the Cauchy integral formula only after 
several tens of pages. Does it mean that a proof of this formula has to be 50 pages long?\,---\,see M.~Klazar \cite{klaz}. Thus we will not discuss derivations of the Stirling formula based on the third and fourth expression. 
Speaking of the $(C)\int$, it is often defined by reduction to the $(R)\int$ or $(RS)\int$ for the real and imaginary parts. But it seems sensible (integration contours are 
usually composed only of straight segments and circular arcs) to integrate these parts just by the (generalized) $(N)\int$, as it is done for example in the textbook \cite[Kapitola 1.6]{vese_komp} of
J. Vesel\'y. 

The simplest integral sufficient for our task is the historically first integral, the $(N)\int$ of I. Newton. It is not a big surprise because in practice  we compute most $(R)\int$s and $(RS)\int$s by the $(N)\int$. Our contribution to the debate (see, for example, B.\,S. Thomson \cite{thom,thom1}) about the merits of the $(R)\int$, the $(RS)\int$, the $(L)\int$ or the 
Henstock--Kurzweil $(HK)\int$ is that the primordial $(N)\int$ is in its way superior because it completely suffices without any further sophistication 
for deducing the fundamental Stirling formula. We develop all properties of the $(N)\int$ needed for these deductions in Section 2. The 
simplest version of the $(N)\int$ for continuous functions suffices for our purposes, for a more general $(N)\int$ with generalized primitives see J. Vesel\'y \cite{vese}
or B.\,S. Thomson \cite{thom,thom1}. 

The two derivations of the Stirling formula in Section 3 are well known to researchers in asymptotic analysis, and so are the results on the 
$(N)\int$ in Section 2 to real analysts, except possibly for Theorems~\ref{fubi_infi_inter} and \ref{fubi_genthm} which are Fubini-type results for iterated Newton integrals over infinite intervals.
Hopefully their combination, presented here with all details, together with the fact that the $(N)\int$ suffices for these derivations may be of interest to both groups and may constitute our original contribution to the subject. 
We want to present all relevant details and are inspired in this by formalized 
mathematics, see \cite[90. Stirling's formula]{formal} for formalizations of the Stirling formula. For example, the Coq formalization builds on the $(R)\int$.
Because of the space and effort limitations we also take some things for granted. It includes the following basic results from real analysis: the definition and properties of the real numbers $\R$, 
the properties of derivatives (the Leibniz formula, differentiation of composite and inverse functions), Lagrange's mean value theorem, uniform continuity of continuous functions on compact sets, and especially the definitions and properties of the functions $\log x$, $e^x$ and
$\cos x$, and of the number $\pi$. The Stirling formula is a popular topic, and many proofs and derivations can be found in the literature, most of them using integrals. We mention a sample of ten: A.\,J.~Coleman \cite{cole}, P. Diaconis and D.~Freedman
\cite{diac_free}, C. Impens \cite{impe}, G.\,J.\,O.~Jameson \cite{jame}, H. Lou \cite{lou}, R. Michel \cite{mich}, M.\,R. Murty and K. Sampath \cite{murt_samp}, S. Niizeki 
and M.~Araki \cite{niiz_arak}, J.\,M. Patin \cite{pati}, and T. Tao \cite{tao}. This list could be much extended.

\section{The Newton integral}\label{sect_intr_newt}

We use the extended reals $\R^*=\R\cup\{-\infty,+\infty\}$ where $\R$ are the real numbers and $-\infty<a<+\infty$ for every $a\in\R$. By an {\em interval $I$} we mean any subset $I\sus\R$ containing more than one element
and such that $a\le x\le b$ with $a,b\in I$ and $x\in\R$ implies $x\in I$. For $a,b\in\R^*$ with $a<b$ we write $(a,b)=\{x\in\R\;|\;a<x<b\}$ for the {\em open intervals}. The {\em compact intervals} are $[a,b]=\{x\in\R\;|\;a\le x\le b\}$ for $a,b\in\R$, $a<b$. Recall that if $I$ is an interval, $F,f\cc I\to\R$ are two functions, and 
$F'(x)=f(x)$ for every $x\in I$, where $F'(x)$ means the corresponding one-sided derivative of $F$ if $x$ is an endpoint of $I$, then $F$ is called a {\em primitive to $f$ (on $I$)}.
For the discussion of complexity of finding or recovering primitives see the studies \cite{doug_kech} by R. Dougherty and A.\,S. Kechris and \cite{frei} by Ch. Freiling, or the surveys \cite{bull,bull1} by P.\,S. Bullen.

\begin{defi}[the Newton integral]\tec\label{defi_newt}
Suppose that $a,b\in\R^*$, $a<b$, and that  $f\cc(a,b)\to\R$ is a real function. The Newton integral of $f$ over 
$(a,b)$ is the real number
$$
(N)\int_a^bf:=F(b^-)-F(a^+)
$$
where $F$ is on $(a,b)$ primitive to $f$ and both limits $F(a^+):=\lim_{x\to a^+}F(x)$ and $F(b^-):=\lim_{x\to b^-}F(x)$ are finite. We set
$$
(N)\int_b^af:=-\;(N)\int_a^bf\;.
$$
\end{defi}
If $f$ does not have a primitive on $(a,b)$ or one of the limits of $F$ does not exist or is infinite, the Newton integral of $f$ is undefined. It is well known and easy to prove by Lagrange's 
mean value theorem that any two primitives to the same function only differ by a constant shift, and thus the definition is correct (independent of the choice of $F$). The functions $F$ and $f$ may be defined also outside $(a,b)$, and therefore the limits of $F$ have to be marked as one-sided. We use the traditional notation
$$
\int_a^b f\;dx=\int_a^b f(x)\;dx
$$
only in situations when it is necessary to identify the integration variable ($x$ in this case). 

Existence of the $(N)\int_a^bf$ for any finite $a,b$ and any $f$ continuous on $[a,b]$ follows from the next theorem. Recall that a sequence of functions $f_n$, $n\in\N$, defined on a set $M\sus\R$ {\em converges on $M$ locally uniformly to a function
$f\cc M\to\R$}, briefly written {\em locally $f_n\rightrightarrows f$ on $M$}, if for every $a\in M$ there is an open interval $I\ni a$ such that for every $\ep>0$ there is an $n_0\in\N$ such that if $n\ge n_0$ and $x\in I\cap M$
then $|f_n(x)-f(x)|<\ep$. If one may always set $I=\R$, we say that {\em $f_n$ converge on $M$ uniformly to $f$} and write briefly {\em $f_n\rightrightarrows f$ on $M$}.

\begin{thm}[primitive by limit transition]\tec\label{thm_prim_limi}
Let $I$ be an interval, $a$ in $I$ be arbitrary but fixed, and functions $f,f_n\cc I\to\R$, $n\in\N$, be such that (i) locally $f_n\rightrightarrows f$ on $I$ and (ii) each $f_n$ has on $I$ a primitive. Then the
primitives $F_n$ to $f_n$ satisfying $F_n(a)=0$, $n\in\N$, converge on $I$ locally uniformly to a primitive $F$ to $f$.
\end{thm}
\proof
First we show that the sequence $F_n(x)$, $n=1,2,\ds$, is Cauchy, uniformly in $x\in J$ for any compact interval $J\sus I$ containing $a$. Indeed, if $m\ge n$ and $x\in J$ then, by Lagrange's mean value theorem,
\begin{eqnarray*}
|F_m(x)-F_n(x)|&\le&|(F_m-F_n)(x)-(F_m-F_n)(a)|+|F_m(a)-F_n(a)|\\
&=&|(x-a)\cdot(f_m-f_n)(b)|\;,
\end{eqnarray*}
for some $b$ lying between $x$ and $a$ and thus in $J$. By (i) and the compactness of $J$,  $f_n\rightrightarrows f$ on $J$ and the last absolute value is for large $n$ uniformly small. Thus for any $\ep>0$ 
there is an $n_0\in\N$ such that if $m\ge n\ge n_0$ then $|F_m(x)-F_n(x)|<\ep$ for every $x\in J$. It follows that for some $F\cc I\to\R$ we have $F_n\rightrightarrows F$ on $J$, and hence locally $F_n\rightrightarrows F$ on $I$.

Next we show that $F$ is on $I$ primitive to $f$. Let an $x_0\in I$ be given and let $J\sus I$ be a compact interval containing $x_0$ in its relative interior. Let an $\ep>0$ be given. Since $f_n\rightrightarrows f$ on $J$, we can take an 
$n_0\in\N$ such that if $m\ge n\ge n_0$ then $|f_m(x)-f_n(x)|<\ep$ for every $x\in J$. We fix an $n\ge n_0$ such that $|f_n(x_0)-f(x_0)|<\ep$. Since $F_n'=f_n$ on $I$, we can take a relatively open interval $K\sus J$ such that $x_0\in K$ 
and for every $x\in K$, $x\ne x_0$, we have $|\frac{F_n(x)-F_n(x_0)}{x-x_0}-f_n(x_0)|<\ep$. Let an $x\in K$, $x\ne x_0$, be given. We fix an $m\ge n$ such that
$\left|\frac{F(x)-F(x_0)}{x-x_0}-\frac{F_m(x)-F_m(x_0)}{x-x_0}\right|<\ep$.
Then for the given $x\in K$ we have, by the previous choices, by Lagrange's mean value theorem, and by the triangle inequality, 
\begin{eqnarray*}
&&\left|\frac{F(x)-F(x_0)}{x-x_0}-f(x_0)\right|\le\left|\frac{F(x)-F(x_0)}{x-x_0}-\frac{F_m(x)-F_m(x_0)}{x-x_0}\right|+\\
&&+\;\left|\frac{(F_m-F_n)(x)-(F_m-F_n)(x_0)}{x-x_0}\right|+\left|\frac{F_n(x)-F_n(x_0)}{x-x_0}-f_n(x_0)\right|+\\
&&+\;|f_n(x_0)-f(x_0)|\\
&&<\ep+|(f_m-f_n)(y)|+\ep+\ep<4\ep\;,
\end{eqnarray*}
for some $y$ lying between $x_0$ and $x$ and thus in $J$. Hence $F'(x_0)=f(x_0)$.
\eproof

\noindent
The following is an existence theorem for the $(N)\int$ on which we rely in the case of bounded intervals.

\begin{cor}[existence of the $(N)\int$]\tec\label{prop_ex_newt}
Let $a,b\in\R$, $a<b$, and $f\cc[a,b]\to\R$ be a continuous function. Then $f$ has a primitive $F$ on $[a,b]$ and the $(N)\int_a^bf$ exists. 
\end{cor}
\proof
It suffices to prove the existence of $F$ because $F(a^+)=F(a)$ and $F(b^-)=F(b)$ by its continuity. Due to compactness of $[a,b]$ the function $f$ is uniformly continuous. 
So for every $n\in\N$ there is a partition $a=a_0<a_1<\ds<a_k=b$ of $[a,b]$ (we do not mark the dependence on $n$) such that 
$$
a_i\le x\le a_{i+1}\Rightarrow|f(x)-f(a_i)|<\frac{1}{n}\;\mbox{ and }\;|f(x)-f(a_{i+1})|<\frac{1}{n}
$$
for every $i=0,1,\ds,k-1$. Let $f_n\cc[a,b]\to\R$ be the piecewise linear continuous function whose graph is the broken line with breaks exactly in the points $(a_i,f(a_i))$, $i=0,1,\ds,k$. 
We check that $f_n$ and $f$ satisfy both hypotheses in Theorem~\ref{thm_prim_limi}. By the definition of $f_n$, if $x\in[a_i,a_{i+1}]$ then the value $f_n(x)$ lies between 
$f(a_i)$ and $f(a_{i+1})$, thus $|f(x)-f_n(x)|<\frac{2}{n}$ and we see that even $f_n\rightrightarrows f$ on $[a,b]$ and (i) holds. Since
for every $u,v,w\in\R$ the function $(u/2)x^2+vx+w$ is primitive on any interval to the linear function $ux+v$, it is easy by employing the shifts $w$ to patch from the local primitives to the linear pieces of $f_n$ on 
the intervals $[a_i,a_{i+1}]$, $i=0,1,\ds,k-1$, a function $g_n$ that is primitive to $f_n$ on the whole interval $[a,b]$. In this we use the fact that for any real function
$h$, if $h'_-(x)=y$ and $h'_+(x)=y$ then $h'(x)=y$. Thus (ii) holds. By Theorem~\ref{thm_prim_limi}, $f$ has on $[a,b]$ a primitive function $F$.
\eproof

\noindent
As is well known one can obtain a primitive $F$ to $f$ also as the Riemann integral $F(x)=(R)\int_a^x f$, but this goes against the spirit of our article. 
Similar limit constructions of primitives appear, for example, in J. Jost \cite[Chapter 6]{jost} or B.\,S. Thomson \cite[Chapter 1.2]{thom}. If one is interested only in proving
the existence of a primitive to any continuous $f$, then the proof in \cite[Chapter 1.2]{thom} is simpler compared to our argument Theorem~\ref{thm_prim_limi} $\to$ Corollary~\ref{prop_ex_newt}.
By a historical note in \cite[Chapter 1.2]{thom} quoting F.\,A. Medvedev \cite[p. 66]{medv}, it was only in 1905 when H.~Lebesgue provided in \cite{lebe} a $(R)\int$-free 
construction of primitives to continuous functions; up to then some arguments justifying their existence were logically circular, as they obtained a primitive in terms of
an integral that they had earlier defined in terms of a primitive. 

\begin{prop}[Hake's theorem]\tec\label{Hake_thm}
Let $a,b\in\R^*$, $a<b$, and $f\cc(a,b)\to\R$ be a function. Then in 
$$
(N)\int_a^b f=\lim_{c\to b^-}(N)\int_a^c f
$$
if one side is defined and finite, so is the other side and the equality holds. Similar result holds for the limit with $c\to a^+$. 
\end{prop}
\proof
If the left side is defined and finite, it is $F(b^-)-F(a^+)$ where $F$ is on $(a,b)$ primitive to $f$. For any $c\in(a,b)$ then, for the restricted $f$ and $F$, the
$(N)\int_a^c f$ exists and equals $F(c^-)-F(a^+)=F(c)-F(a^+)$ by the continuity of $F$ at $c$. The limit transition $c\to b^-$ then shows that the right
side equals $F(b^-)-F(a^+)$ too.

If the right side is defined and finite, for every $c\in(a,b)$ we have on $(a,c)$ a primitive $F_c$ to to the restricted $f$,  and $(N)\int_a^c f=F_c(c^-)-F_c(a^+)$.
By the property of primitives, we can take such $F_c$ that $F_c(a^+)=0$ for every $c$ in $(a,b)$. Then $a<c<d<b\Rightarrow F_c\sus F_d$ (i.e. $F_d$ extends $F_c$)
and $F=\bigcup_{c\in(a,b)}F_c$ is a primitive to $f$ on $(a,b)$. Then 
\begin{eqnarray*}
\lim_{c\to b^-}(N)\int_a^c f&=&\lim_{c\to b^-}(F_c(c^-)-F_c(a^+))=\lim_{c\to b^-}(F(c)-F(a^+))\\
&=&\lim_{c\to b^-}F(c)-F(a^+)=F(b^-)-F(a^+)\\
&=&(N)\int_a^b f\;.
\end{eqnarray*}
\eproof

\noindent
Since $(N)\int_a^c f$ is not defined for $c>b$, we could write $\lim_{c\to b}$ in the statement.

\begin{prop}[linearity and additivity]\tec\label{prop_linear}
Suppose that $a,b\in\R^*$, $a<b$, and that the integrals $(N)\int_a^bf$ and $(N)\int_a^bg$ exist. Then the following holds. 
\begin{enumerate}
    \item For every $\al,\be\in\R$ the function $h=\al f+\be g$ has Newton integral over $(a,b)$ and 
    $$
    (N)\int_a^bh=\al\cdot(N)\int_a^bf+\be\cdot(N)\int_a^bg\;.
    $$
    \item  For every $c\in(a,b)$ the integrals $(N)\int_a^cf$ and $(N)\int_c^bf$ exist and 
    $$
    (N)\int_a^bf=(N)\int_a^cf+(N)\int_c^bf\;.
    $$
\end{enumerate}
\end{prop}
\proof
1. This follows from the fact that if $F$ and $G$ are on $(a,b)$ primitive to $f$ and $g$, respectively, then $\al F+\be G$ is on $(a,b)$ primitive to $\al f+\be g$, and from linearity of functional limits at $a^+$ and $b^-$.

2. If $F$ is on $(a,b)$ primitive to $f$, it (its restriction) is primitive to (the restricted) $f$ also on $(a,c)$ and on $(c,b)$. The integrals $(N)\int_a^cf$ and $(N)\int_c^bf$ exist  because 
$$
\lim_{x\to c^+}F(x)=\lim_{x\to c^-}F(x)=F(c)
$$
by the continuity of $F$ at $c$. Also, $F(b^-)-F(a^+)=(F(b^-)-F(c^+))+(F(c^-)-F(a^+))$ gives the stated equality.
\eproof

\noindent
 Manipulations of integrals very often use part 1, and we will not always acknowledge it.

\begin{prop}[monotonicity]\tec\label{prop_mono_newt}
Suppose that $a,b\in\R^*$, $a<b$, the integrals $(N)\int_a^bf$ and $(N)\int_a^bg$ exist, and that 
$f(x)\le g(x)$ for every $x\in(a,b)$. Then
$$
(N)\int_a^bf\le(N)\int_a^bg\;.
$$
\end{prop}
\proof
Let $a<a'<b'<b$ where $a',b'\in\R$ and let $F$ and $G$ be on $(a,b)$ primitive to $f$ and $g$, respectively. By Lagrange's mean value theorem we have, with some $c\in(a',b')$,  
\begin{eqnarray*}
G(b')-G(a')-(F(b')-F(a'))&=&(G-F)(b')-(G-F)(a')\\
&=&(b'-a')\cdot(G-F)'(c)\\
&=&(b'-a')(g(c)-f(c))\\&\ge&0
\end{eqnarray*} 
because $g-f\ge0$ on $(a,b)$. So 
$$
F(b')-F(a')\le G(b')-G(a')\;,
$$
and limit transitions $a'\to a^+$ and $b'\to b^-$ give the stated inequality. 
\eproof

\noindent
As a corollary we obtain the most often used estimate in the integral calculus.

\begin{cor}[ML bound]\tec\label{ML_bound}
Suppose that $a,b,c\in\R$, $a<b$, the integral $(N)\int_a^bf$ exists, and $f(x)\le c$ (resp. $f(x)\ge c$) for every $x\in(a,b)$. Then
$$
(N)\int_a^bf\le c(b-a)\ \bigg(\mbox{resp. }(N)\int_a^bf\ge c(b-a)\bigg)\;.
$$
\end{cor}
\proof
Apply the proposition to $f(x)$ and the constant function $c$, and compute the $(N)\int$ of a constant function.
\eproof

The next theorem can be found in a more general form with generalized primitives in J. Vesel\'y \cite[V\v{e}ta 11.3.13]{vese}.

\begin{thm}[integration by parts]\tec\label{prop_pp_newt}
Let $a,b\in\R^*$, $a<b$, and $F$, resp. $G$, be on $(a,b)$ primitive to $f$, resp. to $g$. Then in
$$
(N)\int_a^b fG=\big((FG)(b^-)-(FG)(a^+)\big)-(N)\int_a^b Fg
$$
if two of the three terms are defined and finite then so is the third one and the equality holds.
\end{thm}
\proof
Suppose that the first term $E(b^-)-E(a^+)$, where $E$ is on $(a,b)$ primitive to $fG$, and the second term $(FG)(b^-)-(FG)(a^+)$ are defined and finite. By the Leibniz rule,
$$
(FG-E)'=fG+Fg-fG=Fg
$$
on $(a,b)$ and $FG-E$ is primitive to $Fg$. Also, by the assumptions, $(FG-E)(b^-)=(FG)(b^-)-E(b^-)$ and $(FG-E)(a^+)=(FG)(a^+)-E(a^+)$. The stated equality therefore follows by subtraction and rearrangement. If the third term and the second 
term are defined and finite, the argument is similar. Suppose that the first term $E(b^-)-E(a^+)$ and the third term $D(b^-)-D(a^+)$, where $E$ and $D$ are on $(a,b)$ primitive to $fG$ and $Fg$, respectively, are defined 
and finite. By the Leibniz rule,
$$
(E+D)'=fG+Fg=(FG)'
$$
on $(a,b)$. Thus (by Lagrange's mean value theorem) $E+D$ and $FG$ only differ by a constant shift $c$. Hence $(FG)(b^-)=E(b^-)+D(b^-)+c$ and $(FG)(a^+)=E(a^+)+D(a^+)+c$. The stated equality again follows by subtraction 
and rearrangement.
\eproof

\begin{prop}[substitution rule]\tec\label{prop_subst}
Suppose that $a,b,c,d\in\R^*$, $a<b$ and $c<d$, $g\cc (c,d)\to(a,b)$, $g(x)\to a$ for $x\to c$, $g(x)\to b$ for $x\to d$, $g$ is differentiable on $(c,d)$, $f\cc(a,b)\to\R$, and the 
$(N)\int_a^bf$ exists. Then the next integral exists and
$$
(N)\int_c^d (f\circ g)g'=(N)\int_a^b f=(N)\int_{g(c)}^{g(d)} f\;.
$$
We extend $g$ by $g(c)=a$ and $g(d)=b$ by limit transitions, and understand it as a mere notation when $c=-\infty$ or $d=+\infty$.
\end{prop}
\proof
Let $F$ be on $(a,b)$ primitive to $f$. Then 
$$
(N)\int_a^b f=F(b^-)-F(a^+)=(F\circ g)(d^-)-(F\circ g)(c^+)=(N)\int_c^d (f\circ g)g'
$$
because $(F\circ g)'=(f\circ g)g'$ on $(c,d)$ and therefore $F\circ g$ is on $(c,d)$ primitive to $(f\circ g)g'$. 
\eproof

\noindent
In the situation when the substitution $g$ flips the interval by $g(x)\to b$ for $x\to c$  and $g(x)\to a$ for $x\to d$ and we modify the hypothesis accordingly, we obtain identical formula:
\begin{eqnarray*}
(N)\int_{g(c)}^{g(d)} f&=&(N)\int_b^a f=-\;(N)\int_a^b f
=-F(b^-)+F(a^+)\\
&=&-\;(F\circ g)(c^+)+(F\circ g)(d^-)=(N)\int_c^d (f\circ g)g'\;.
\end{eqnarray*}

For the second derivation of the Stirling formula we need for the Newton integral a Fubini-type result. We obtain it in the next two theorems. The first one appears, with a different proof, in P. Walker \cite[Theorem A.9 (i)]{walk} 
for the $(W)\int$. The integral, which we call tentatively Walker's, is introduced in \cite[Chapter 4]{walk} that was not available to us, and we could not determine its relation to other integrals. Later we see that $(W)\int\ne(N)\int$.

\begin{thm}[Fubini \`a la Newton]\tec\label{thm_fubi}
Suppose that $a,b,c,d\in\R$, $a<b$ and $c<d$, and
$$
f=f(x,\,y)\cc[a,\,b]\times[c,\,d]\to\R
$$ 
is a continuous function. Then the following two iterated Newton integrals exist and are equal:
$$
(N)\int_a^b\left((N)\int_c^d f(x,\,y)\,dy\right)\,dx=(N)\int_c^d\left((N)\int_a^b f(x,\,y)\,dx\right)\,dy\;.
$$
\end{thm}
\proof
Each inner integral $I(x)=(N)\int_c^d f(x,y)\,dy$ exists by Corollary~\ref{prop_ex_newt}. If $x_1,x_2\in[a,b]$ then
$$
I(x_1)-I(x_2)=(N)\int_c^d(f(x_1,\,y)-f(x_2,\,y))\,dy\;.
$$
By the uniform continuity of $f(x,y)$ on the compact rectangle $[a,\,b]\times[c,\,d]$ we see that for close $x_1$ and $x_2$ the value
$|f(x_1,y)-f(x_2,y)|$ is small for any $y$, and thus by Corollary~\ref{ML_bound} also the integral and $|I(x_1)-I(x_2)|$ are small\,---\,$I(x)$ is continuous. Thus $(N)\int_a^b I(x)$ exists. Similar argument shows existence of the integrals 
$J(y)=(N)\int_a^bf(x,y)\,dx$ and $(N)\int_c^d J(y)$ on the right side of the formula. 

We prove the equality by showing that the two iterated Newton integrals are arbitrarily close. Let $\ep>0$ be given. By the uniform continuity of $f$ on the rectangle there exist a partition $a=a_0<a_1<\ds<a_k=b$ of $[a,b]$, a partition $c=b_0<b_1<\ds<b_l=d$ of $[c,d]$,
and constants $c_{i,j}\in\R$, $i=0,1,\ds,k-1$ and $j=0,1,\ds,l-1$, such that if $(x,y)$ lies in $[a_i,a_{i+1}]\times[b_j,b_{j+1}]$ then $|f(x,y)-c_{i,j}|<\ep$. Let $f_{i,j}$ be the restriction of $f$ to $[a_i,a_{i+1}]\times[b_j,b_{j+1}]$ and let $c_{i,j}$ 
also denote the constant function $c_{i,j}$ on this rectangle. Then
\begin{eqnarray*}
I_{i,j}&:=&(N)\int_{a_i}^{a_{i+1}}\left((N)\int_{b_j}^{b_{j+1}}c_{i,j}\,dy\right)\,dx=(N)\int_{a_i}^{a_{i+1}}c_{i,j}(b_{j+1}-b_j)\,dx\\
&=&c_{i,j}(a_{i+1}-a_i)(b_{j+1}-b_j)=(N)\int_{b_j}^{b_{j+1}}c_{i,j}(a_{i+1}-a_i)\,dy\\
&=&(N)\int_{b_j}^{b_{j+1}}\left((N)\int_{a_i}^{a_{i+1}}c_{i,j}\,dx\right)\,dy=:J_{i,j}\;.
\end{eqnarray*}
By parts 1 and 2 of Proposition~\ref{prop_linear}, 
\begin{eqnarray*}
(N)\int_a^b\left((N)\int_c^d f\,dy\right)\,dx&=&(N)\int_a^b\left(\sum_{j=0}^{l-1}(N)\int_{b_j}^{b_{j+1}} f\,dy\right)\,dx\\
&=&\sum_{j=0}^{l-1}(N)\int_a^b\left((N)\int_{b_j}^{b_{j+1}} f\,dy\right)\,dx\\
&=&\sum_{j=0}^{l-1}\sum_{i=0}^{k-1}(N)\int_{a_i}^{a_{i+1}}\left((N)\int_{b_j}^{b_{j+1}} f_{i,j}\,dy\right)\,dx\;.
\end{eqnarray*}
A similar computation shows that 
$$
(N)\int_c^d\left((N)\int_a^b f\,dx\right)\,dy=\sum_{i=0}^{k-1}\sum_{j=0}^{l-1}(N)\int_{b_j}^{b_{j+1}}\left((N)\int_{a_i}^{a_{i+1}} f_{i,j}\,dx\right)\,dy\;.
$$
Since $|f_{i,j}-c_{i,j}|<\ep$ on $[a_i,a_{i+1}]\times[b_j,b_{j+1}]$, it follows by Corollary~\ref{ML_bound} that the first iterated Newton integral in the equality we are proving differs from $\sum_{j=0}^{l-1}\sum_{i=0}^{k-1}I_{i,j}$ 
by less than $\ep(b-a)(d-c)$, and the second one  differs from $\sum_{i=0}^{k-1}\sum_{j=0}^{l-1}J_{i,j}$ by less than $\ep(d-c)(b-a)$. Since always $I_{i,j}=J_{i,j}$, these two double sums are equal and the two iterated Newton integrals differ by less
than  $2\ep(b-a)(d-c)$, as we need.
\eproof

\noindent
The theorem inverts the well known result that $\partial_x\partial_yf=\partial_y\partial_xf$ at a point if both second order partial derivatives exist in a neighborhood of 
the point and are continuous at it. 

But what we need in Section 3 is an extension of the previous theorem with infinite $b$ and $d$. Then \cite[Theorem A.9 (ii)]{walk} says:
\begin{quote}
    (ii) If $f$ is continuous on $I\times J$ where $I,J$ are intervals in $R$ which may be finite or infinite, and if $f$ is positive on $E$ 
    [$=I\times J$] then the integrals in (A.1) [the two iterated integrals] are either all infinite, or all finite and equal.
\end{quote}
This does not hold for the $(N)\int$. Consider a continuous function 
$$
f=f(x,\,y)\cc[0,\,+\infty)^2\to(0,\,1]
$$ 
such that $f(x,1)=1$ for every $x\ge0$, and such that for 
each fixed $x\ge0$ the section $f(x,y)$ first increases for $0\le y\le 1$ from $0^+$ to $1$ and then for $y\ge 1$ rapidly decreases 
from $1$ to $0^+$ in such a way that the width of the base of this peak decreases for $x\to+\infty$ to $0$ fast enough so that each $J(x)=(N)\int_0^{+\infty}f(x,y)\,dy$ exists, $J\cc[0,+\infty)\to(0,+\infty)$ is 
continuous, and
$$
(N)\int_0^{+\infty}J(x)
$$
exists. Then the iterated Newton integral
$$
(N)\int_0^{+\infty}\left((N)\int_0^{+\infty} f(x,\,y)\,dy\right)\,dx=(N)\int_0^{+\infty}J(x)
$$
exists. However, the other iterated Newton integral
$$
(N)\int_0^{+\infty}\left((N)\int_0^{+\infty} f(x,\,y)\,dx\right)\,dy
$$
is undefined because for $y=1$ the inner Newton integral is not defined. The reader will have no problems to supply numerical details. 

We do not give a general Fubini-type theorem for the $(N)\int$ over infinite intervals strong enough to prove 
Proposition~\ref{gauss_int} because we could not find such a theorem. Instead we directly establish only the needed instance 
for the function $ue^{-u^2(1+v^2)}$.

\begin{thm}[a (N) Fubini result]\tec\label{fubi_infi_inter}
The next two iterated Newton integrals exist and are equal: if $f(x,z)=xe^{-x^2(1+z^2)}=xe^{-x^2}e^{-x^2z^2}$ then
$$
(N)\int_0^{+\infty}\left((N)\int_0^{+\infty}f(x,\,z)\,dz\right)\,dx=(N)\int_0^{+\infty}\left((N)\int_0^{+\infty}f(x,\,z)\,dx\right)\,dz\;.
$$
\end{thm}
\proof
The $(N)\int_0^{+\infty} e^{-x^2}=:c>0$ exists by the majorization $e^{-x^2}\le e^{-x}$ for $x\ge1$, Corollary~\ref{prop_ex_newt}, and Propositions~\ref{Hake_thm}, \ref{prop_linear} (part 2), and \ref{prop_mono_newt}. The calculation
\begin{eqnarray*}
(N)\int_0^{+\infty} e^{-x^2}\cdot(N)\int_0^{+\infty} e^{-y^2}&=&(N)\int_0^{+\infty} e^{-x^2}\left((N)\int_0^{+\infty} e^{-y^2}\right)\,dx\\
&=&(N)\int_0^{+\infty}\left((N)\int_0^{+\infty} e^{-x^2-y^2}\,dy\right)\,dx\\
&=&(N)\int_0^{+\infty}\left((N)\int_0^{+\infty}f(x,\,z)\,dz\right)\,dx
\end{eqnarray*}
then shows that the first iterated Newton integral $A$ exists. On the first two lines we multiplied an integral by a constant
according to part 1 of Proposition~\ref{prop_linear} and we used that $e^a e^b=e^{a+b}$. On the third line we used Proposition~\ref{prop_subst} with the substitution $y\leftarrow z$, $y=xz$. To prove the equality we estimate how much the last iterated integral $A$
differs from its finite approximation
$$
A(b):=(N)\int_0^b\left((N)\int_0^b f(x,\,z)\,dz\right)\,dx,\ \R\ni b\ge1\;.
$$
The integrals $A(b)$ exist by Theorem~\ref{thm_fubi}.

For any $b\ge1$ (we justify the estimates after the computation),
\begin{eqnarray*}
&&0\le A-A(b)=(N)\int_0^b\left((N)\int_b^{+\infty}f(x,z)\,dz\right)\,dx\,+\\
&&+\,(N)\int_b^{+\infty}\left((N)\int_0^{+\infty}f(x,z)\,dz\right)\,dx\\
&&\le(N)\int_0^b\left((N)\int_{bx}^{+\infty}e^{-x^2}e^{-y^2}\,dy\right)\,dx+(N)\int_b^{+\infty}\left((N)\int_0^1 c_0x^{-2}\,dz\,+\right.\\
&&\left.+\,(N)\int_1^{+\infty} c_0x^{-2}e^{-z}\,dz\right)\,dx\\
&&\le(N)\int_0^{1/b^{1/2}} c\,dx+(N)\int_{1/b^{1/2}}^b\left((N)\int_{b^{1/2}}^{+\infty}e^{-y}\right)\,dx+\\
&&+\,(N)\int_b^{+\infty}(c_0x^{-2}+c_0x^{-2}e^{-1})\,dx\\
&&\le cb^{-1/2}+be^{-b^{1/2}}+c_0(1+e^{-1})b^{-1}\;.
\end{eqnarray*}
In the initial $=$ we used part 2 of Proposition~\ref{prop_linear}. In the next $\le$ we returned from $z$ to the variable $y$ and set $c_0=\max_{x\ge 1}xe^{-x^2}/x^{-2}$. In the penultimate $\le$ we invoked the existence of $(N)\int_0^{+\infty}e^{-y^2}$. We also were using part 2 of 
Proposition~\ref{prop_linear}, Proposition~\ref{prop_mono_newt}, Definition~\ref{defi_newt}, and majorizations $e^{-a}\le1$ for $a\ge0$ and $e^{-a^2}\le e^{-a}$ for $a\ge1$. Thus $A-A(b)\to0$ as $b\to+\infty$.

By Theorem~\ref{thm_fubi},
$$
A(b)=B(b):=(N)\int_0^b\left((N)\int_0^b f(x,\,z)\,dx\right)\,dz
$$
for any $b\in\R$ with $b>0$. We complete the proof by showing that 
$$
B:=(N)\int_0^{+\infty}\left((N)\int_0^{+\infty}f(x,\,z)\,dx\right)\,dz
$$
exists and that $B-B(b)\to0$ as $b\to+\infty$. 

For any $z,b>0$ we define $I(z,b)=(N)\int_0^b f(x,z)\,dx$, this integral exists by Corollary~\ref{prop_ex_newt}. Since $f(x,z)$ is for $x\ge1$ majorized by $xe^{-x}$, the inner integral
$$
I(z):=(N)\int_0^{+\infty}f(x,\,z)\,dx=\lim_{b\to+\infty}I(z,\,b)
$$
exists for any $z\ge0$ by Propositions~\ref{prop_linear} (part 2), \ref{prop_mono_newt}, and \ref{Hake_thm}. We prove that $I(z)$ is continuous for
$z\ge0$. By the uniform continuity of $f(x,z)$ on compact sets, for any given $z_0\ge0$, $b\ge1$, and $\ep>0$ there is a $\de>0$ such that if $z\ge0$ satisfies $|z-z_0|<\de$ then
$|f(x,z)-f(x,z_0)|<\ep$ for any $x\in[0,b]$. Then
$$
|I(z)-I(z_0)|<(N)\int_0^b\ep\,dx+(N)\int_b^{+\infty}xe^{-x}=b\ep+be^{-b}+e^{-b}\;,
$$
which shows that $I(z)$ is continuous at $z_0$. Thus $(N)\int_0^bI(z)$ exists for any $b>0$ by Corollary~\ref{prop_ex_newt}.
Let $c_1=\max_{x\ge0}xe^{-x^2}>0$. For any $z\ge1$ we have
\begin{eqnarray*}
0\le I(z)&<&(N)\int_0^{1/z^{2/3}}x\cdot1\,dx+(N)\int_{1/z^{2/3}}^{+\infty}xe^{-x^2}e^{-xz}\,dx\\
&<&z^{-4/3}+c_1e^{-z^{1/3}}/z<c_2z^{-4/3}\;,
\end{eqnarray*}
for an absolute constant $c_2$. This majorizations implies that the integral
$$
B=(N)\int_0^{+\infty}I(z)
$$
exists. 

It remains to estimate its distance from $B(b)$. For any $b\ge1$ we have, using again part 2 of Proposition~\ref{prop_linear}, that
\begin{eqnarray*}
&&0\le B-B(b)=(N)\int_0^b\left((N)\int_b^{+\infty}f(x,z)\,dx\right)\,dz+\\
&&+\,(N)\int_b^{+\infty}\left((N)\int_0^{+\infty}f(x,z)\,dx\right)\,dz\\
&&<(N)\int_0^b\left((N)\int_b^{+\infty}xe^{-x}\right)\,dz+(N)\int_b^{+\infty}I(z)\\
&&<(N)\int_0^b(be^{-b}+e^{-b})\,dz+(N)\int_b^{+\infty}c_2z^{-4/3}\\
&&=(b^2+b)e^{-b}+3c_2b^{-1/3}\;,
\end{eqnarray*}
and again $B-B(b)\to0$ as $b\to+\infty$. Thus $A=B$, the two iterated Newton integrals are equal.
\eproof

\noindent
We hoped to prove Proposition~\ref{gauss_int} as an instance of a general Fubini theorem for iterated Newton integrals of $f(x,y)$ over $[0,+\infty)^2$, when $f(x,y)$
satisfies a symmetric decay condition for $x,y\to+\infty$. But we only could prove (again as a corollary of Theorem~\ref{thm_fubi}) the following theorem which unfortunately does not apply to 
 $f(x,y)=xe^{-x^2}e^{-x^2y^2}$; we omit the proof.

\begin{thm}[a (N) Fubini theorem]\tec\label{fubi_genthm}
If $c>0$ is a constant and $f=f(x,y)\cc[0,+\infty)^2\to\R$ is a continuous function such that 
$$
|f(x,y)|\le c\max(x,y)^{-3}\mbox{ for }\max(x,y)\ge1\;,
$$ 
then the next two iterated Newton integrals exist and are equal,
$$
(N)\int_0^{+\infty}\left((N)\int_0^{+\infty}f(x,\,y)\,dy\right)\,dx=(N)\int_0^{+\infty}\left((N)\int_0^{+\infty}f(x,\,y)\,dx\right)\,dy\;.
$$
\end{thm}

\section{The Stirling formula}

With the help of the properties of the $(N)\int$ in Section 2 we prove in two ways the next basic asymptotic formula.
\begin{thm}[Stirling formula]\tec\label{thm_stir}
For $n\in\N$ one has
\begin{eqnarray*}
n!&=&\prod_{i=1}^n i=\#\{(a_1,\,\ds,\,a_n)\in\{1,\,2,\,\ds,\,n\}^n\;|\;|\{a_1,\,\ds,\,a_n\}|=n\}\\
&=&(1+o(1))\sqrt{2\pi n}\left(\frac{n}{e}\right)^n
\end{eqnarray*}
where $\pi=3.14159\ds$ and $e=2.71828\ds$ are well known constants.
\end{thm}
Here the asymptotic notation $f=o(g)$ for $f,g\cc M\to\R$, $M\sus\R$ and $\sup(M)=+\infty$, means that 
$$
\lim_{x\to+\infty}\frac{f(x)}{g(x)}=0\;.
$$

We start the first proof by borrowing an estimate, but not its proof, from G.~Tenenbaum \cite[Theorem I.0.4]{tene}. 
There it is proven via integration by parts in a $(RS)\int$. We actually learned this proof of Theorem~\ref{thm_stir} from 
G.~Tenenbaum \cite[Exercise 3 on p. 8]{tene}. We could do without the next proposition, see the remark on telescoping after Corollary~\ref{suma},
but we keep it as a basic result on the interplay of sums and integrals. $\Z$ denotes the ring of integers.

\begin{prop}[basic estimate]\tec\label{basi_esti}
Let $a,b\in\Z$, $a<b$, and $f\cc[a,b]\to\R$ be a continuous monotonic function. Then there exists a number $\theta\in[0,1]$ such that
$$
\sum_{a<n\le b}f(n)=(N)\int_a^bf+\theta(f(b)-f(a))\;.
$$
\end{prop}
\proof
Suppose that $f$ is nondecreasing, the proof for nonincreasing $f$ is similar (by reverting the next two inequalities). The equality we need to prove is equivalent with the estimate
$$
0\le\sum_{a<n\le b}f(n)-(N)\int_a^bf\le f(b)-f(a)\;.
$$
Note that by part 2 of Proposition~\ref{prop_linear}, the estimate is additive: if $c$ is an integer with $a<c<b$ and we have the estimate for both pairs
$a,c$ and $c,b$ (in place of $a,b$), then by summing we get it for $a,b$. Therefore it suffices  to prove it only for $b=a+1$ (one can partition $[a,b]$ into unit 
intervals $[x,x+1]$, $x\in\Z$ with $a\le x<b$). For $b=a+1$ the estimate becomes
$$
0\le f(a+1)-(N)\int_a^{a+1}f\le f(a+1)-f(a)\;.
$$
By Corollary~\ref{ML_bound},
$$
f(a)\cdot1\le(N)\int_a^{a+1}f\le f(a+1)\cdot1\;,
$$
and the instance $a,a+1$ of the estimate follows.
\eproof

\noindent
The  more general result \cite[Theorem I.0.4]{tene} drops the continuity of $f$ and uses the $(R)\int_a^b f$. Then one can prove it easily by lower and upper Riemann--Darboux sums, which seems to be the simplest of the three arguments (if one has already built the theory of the $(R)\int$).
We learned the additive estimate trick used
in the previous proof in E.\,C. Titchmarsh \cite[p. 13/14]{titc_rf}. By it he gives a simple, few lines proof of the more precise formula ($a,b,c\in\R$ with $a<b$, $f=f(t)\cc[a,b]\to\R$ is continuously differentiable,
and $\{x\}=x-\lfloor x\rfloor\in[0,1)$ is the fractional part of $x\in\R$) 
\begin{eqnarray*}
\sum_{a<n\le b}f(n)&=&(R)\int_a^b f+(R)\int_a^b (\{t\}+c)f'(t)\,dt\\
&&+\,(\{a\}+c)f(a)-(\{b\}+c)f(b)\;.
\end{eqnarray*}
Another proof in a monograph on analytic number theory takes $1\frac{1}{2}$ pages. The Euler--Maclaurin summation formula (EMSF), see for example
\cite[Chapter I.0.2]{tene}, is much more precise. An alternative to EMSF, using only integrals and with derivatives only in the error term, was recently proposed by I. Pinelis \cite{pine}. 

We use the standard asymptotic notation $O$ and $\ll$: if $f,g\cc M\to\R$, $M\sus\R$, then $f=O(g)$ (on $M$) and $f\ll g$ (on $M$) both mean that there is a constant $c>0$ such that
for every $x\in M$ ones has $|f(x)|\le c|g(x)|$.

\begin{cor}[reciprocal squares]\tec\label{suma}
For all $n\in\N$ one has
$$
\sum_{m=1}^n O(m^{-2})=c+O(n^{-1})\;,
$$
for a constant $c\in\R$.
\end{cor}
\proof
The claim is that if $f\cc\N\to\R$ satisfies $f(m)=O(m^{-2})$ then the sum $\sum_{m=1}^n f(m)$ has the stated asymptotics. We have
$$
\sum_{m=1}^n f(m)=\lim_{N\to\infty}\sum_{m=1}^N f(m)-\lim_{N\to\infty}\sum_{m=n+1}^N f(m)=:\lim_{N\to\infty}S(N)-\lim_{N\to\infty}S(n,N)
$$
provided, of course, that both limits exist and are finite. But for $M>N$ we have by the previous proposition that
$$
|S(M)-S(N)|\ll\sum_{m=N+1}^M m^{-2}\le(N)\int_N^M x^{-2}=N^{-1}-M^{-1}<N^{-1}\;.
$$
Thus $S(N)$, $N=1,2,\ds$, is a Cauchy sequence and has a finite limit $c$. Similar argument shows for each $n$ existence and finiteness of the second limit. By the previous proposition we again have ($N>n$)
$$
|S(n,N)|\ll\sum_{m=n+1}^N m^{-2}\le(N)\int_n^Nx^{-2}=n^{-1}-N^{-1}<n^{-1}\;.
$$
Therefore the second limit is $O(1/n)$.
\eproof

\noindent
Alternatively, we can bound finite sums of reciprocal squares without any integral by using telescoping sums with the telescoper $m^{-2}=m^{-1}-(m+1)^{-1}$.
The previous proof is in a way remarkable. Usually one obtains infinite sums (products, integrals, $\ds$) as
limit cases of finite approximations, one of the best known examples being ($|q|<1$)
$$
\sum_{n=0}^{\infty}q^n=\lim_{n\to\infty}(1+q+q^2+\ds+q^n)=\lim_{n\to\infty}\frac{1-q^{n+1}}{1-q}=\frac{1}{1-q}\;.
$$
%\end{eqnarray*}
In contrast, the previous proof reverts this process and expresses a finite sum by two infinite ones. There seems to be no other way to deduce this
asymptotics apparently involving no infinite expression (the infinity, however, hides in ``all $n\in\N$'') than via the limits at infinity.

In the following proposition we use the Taylor expansion $\log(1+x)=x-\frac{x^2}{2}+O(x^3)$ ($x\in[-\frac{1}{2},2]$) which is yet another application of Lagrange's
mean value theorem. 

\begin{prop}[first expression of $n!$ by an $\int$]\tec\label{appr_sum_logs}
There is a real constant $c$ such that for all $n\in\N$ we have
$$
\log(n!)=c+O(1/n)+(N)\int_{\frac{1}{2}}^{n+\frac{1}{2}}\log x\;.
$$
\end{prop}
\proof
We prove that for all $m\in\N$,
$$
(N)\int_{m-\frac{1}{2}}^{m+\frac{1}{2}}\log x=\log m+O(m^{-2})\;.
$$
Indeed, by Definition~\ref{defi_newt} and by the expansion of $\log(1+x)$ the integral equals
\begin{eqnarray*}
   && (x\log x-x)(m+1/2)-(x\log x-x)(m-1/2)\\
    &&=m\log\left(1+\frac{1}{m-1/2}\right)+\log m+\frac{\log(1-1/4m^2)}{2}-1\\
    &&=\log m+\frac{2m(m-1/2)-m-2(m-1/2)^2}{2(m-1/2)^2}+O(m^{-2})+O(m^{-2})\\
    &&=\log m+\frac{-1/2}{2(m-1/2)^2}+O(m^{-2})=\log m+O(m^{-2})\;.
\end{eqnarray*}
Using equation $\log(n!)=\sum_{m=1}^n\log m$, part 2 of Proposition~\ref{prop_linear} and Corollary~\ref{suma} we get the first expression. 
\eproof

Now the Stirling formula with an undetermined constant follows easily. We use another Taylor expansion $e^x=1+O(x)$ ($x\in[-c,c]$ for any $c>0$) which implies that 
$e^{O(1/n)}=1+O(1/n)$ (for $n\in\N$). 

\begin{prop}[incomplete Stirling formula]\tec\label{stir_unde}
There is a real constant $d>0$ such that for all $n\in\N$ we have
$$
n!=(d+O(1/n))\sqrt{n}\left(\frac{n}{e}\right)^n\;.
$$
\end{prop}
\proof
We compute the integral in the previous proposition in the same way as in its proof and get that
\begin{eqnarray*}
\log(n!)&=&c+O(1/n)+(x\log x-x)(n+1/2)-(x\log x-x)(1/2)\\
&=&n\log(n+1/2)-n+\frac{\log(n+1/2)}{2}+c_0+O(1/n)\\
&=&n\log n-n+\frac{\log n}{2}+n\log(1+1/2n)+\frac{\log(1+1/2n)}{2}+\\
&&+\;c_0+O(1/n)\\
&=&n\log n-n+\frac{\log n}{2}+c_1+O(1/n)\;.
\end{eqnarray*}
We used the above expansion of $\log(1+x)$, collected in the $c_i$ several constant contributions to $c$, and merged several $O(1/n)$ terms in one.
Applying the exponential function we get the expression for $n!$, with $d=e^{c_1}$.
\eproof

\noindent
We remark that if one is in Proposition~\ref{stir_unde} content with $o(1)$ in place of $O(1/n)$, then the argument so far can be shortened and made integral-free by 
simply proving that the sequence $(n!/\sqrt{n}e^{-n}n^n)$ is monotonic and bounded (see for example \cite{impe}). 

It remains to prove that $d=\sqrt{2\pi}$. We do it by another and quite unexpected, at least to the author, application of (Newton) integrals.

\begin{prop}[resolving a recurrence by $\int$]\tec\label{wn}
Suppose that the sequence $(W_n)$ of positive real numbers is given by the recurrence
$$
W_0=\frac{\pi}{2},\ W_1=1,\mbox{ and for $n\ge2$,}\;W_n=\frac{n-1}{n}W_{n-2}\;. 
$$
Then
$$
\lim_{n\to\infty}\frac{W_n}{W_{n-1}}=1\;.
$$
\end{prop}
\proof
The trick is to prove that
$$
W_n=(N)\int_0^{\pi/2}(\cos x)^n,\ n\in\N_0\;.
$$
By Definition~\ref{defi_newt}, $(N)\int_0^{\pi/2}(\cos x)^0=x(\pi/2)-x(0)=\pi/2$ and $(N)\int_0^{\pi/2}\cos x=\sin(\pi/2)-\sin(0)=1$. For $n\ge2$ 
one has by Theorem~\ref{prop_pp_newt} (Corollary~\ref{prop_ex_newt} shows that in the integration by parts identity below both the first and the third term are defined and finite) and by part 1 of Proposition~\ref{prop_linear}
that, denoting $h(x)=(\sin x)(\cos x)^{n-1}$ and using that $\sin^2x=1-\cos^2x$,
\begin{eqnarray*}
&&(N)\int_0^{\pi/2}(\cos x)^n=(N)\int_0^{\pi/2}(\sin x)'(\cos x)^{n-1}\\
&&=h(\pi/2)-h(0)-(N)\int_0^{\pi/2}(\sin x)((\cos x)^{n-1})'\\
&&=0-0+(n-1)\cdot(N)\int_0^{\pi/2}(\sin x)^2(\cos x)^{n-2}\\
&&=(n-1)\left((N)\int_0^{\pi/2}(\cos x)^{n-2}-(N)\int_0^{\pi/2}(\cos x)^n\right)\;.
\end{eqnarray*}
Thus the sequences $\left((N)\int_0^{\pi/2}(\cos x)^n\right)$ and $(W_n)$, $n=0,1,2,\ds$, follow the same recurrence and coincide. Crucially\,---\,this is hard to get 
from the original definition of $W_n$ but it follows easily from the integral representation\,---\,the sequence $(W_n)$ is nonincreasing, it in fact decreases. Indeed, since 
$0\le \cos^n\le\cos^{n-1}$ on $[0,\pi/2]$, Proposition~\ref{prop_mono_newt} shows that $W_n\le W_{n-1}$. Thus for $n\ge2$ we have, by the monotonicity of $W_n$ and the recurrence,
$$
1=\frac{W_{n-1}}{W_{n-1}}\ge\frac{W_n}{W_{n-1}}\ge\frac{W_{n+1}}{W_{n-1}}=\frac{n}{n+1}\;\mbox{ and }\;\frac{W_n}{W_{n-1}}\to1,\ n\to\infty\;.
$$
\eproof

\noindent
Before we complete the determination of $d$ we contemplate for a while the function $\cos x$ used in the previous proof. If to define $\cos x$ 
one needed, say, the $(R)\int$, our undertaking would be less convincing. (We were in a similar situation at the beginning when we needed primitives to continuous functions.)
This function is defined by a limit process but without Riemann integral,
$$
\cos x=\sum_{n=0}^{\infty}\frac{(-1)^nx^{2n}}{(2n)!},\ \ x\in\R\;.
$$
From this formula one derives without using the $(R)\int$ all properties of $\cos x$ needed for the proof, such as the related function $\sin x$, the identity $\sin^2+\cos^2=1$, the relations $\sin'x=\cos x$ and $\cos'x=-\sin x$, and the fact that $\pi/2$ is the smallest positive zero of $\cos x$. Which actually serves
as a definition of $\pi$ for our article. If the adopted definition of $\pi$ were that 
$$
\pi=\lim_{n\to\infty}\frac{n!^2}{2n\left(\frac{n}{e}\right)^{2n}}\;,
$$
we would be done after Proposition~\ref{stir_unde}.

The recurrence for $W_n$ has for $n\in\N$ another explicit solution:
\begin{eqnarray*}
W_{2n}&=&\frac{(2n-1)(2n-3)\ds1}{2n(2n-2)\ds2}\cdot\frac{\pi}{2}=\frac{(2n)!}{(2^n n!)^2}\cdot\frac{\pi}{2}\ \mbox{ and}\\
W_{2n+1}&=&\frac{2n(2n-2)\ds2}{(2n+1)(2n-1)\ds3}\cdot1=\frac{(2^n n!)^2}{(2n+1)!}\;.
\end{eqnarray*}
Employing the asymptotic notation $f\sim g$ which for $f,g\cc M\to\R$, $M\sus\R$ and $\sup(M)=+\infty$, means that 
$$
\lim_{x\to+\infty}\frac{f(x)}{g(x)}=1\;,
$$
we get by Propositions~\ref{wn} and \ref{stir_unde} that
$$
1\sim\frac{W_{2n+1}}{W_{2n}}\sim\frac{(2^n n!)^4}{2n(2n)!^2}\cdot\frac{2}{\pi}\sim\frac{2^{4n}\cdot d^4\cdot n^2\cdot(n/e)^{4n}}{2n\cdot d^2\cdot 2n\cdot(2n/e)^{4n}}\cdot\frac{2}{\pi}=\frac{d^2}{2\pi}\;.
$$
Thus $d=\sqrt{2\pi}$ and the first proof of Theorem~\ref{thm_stir} is complete.\eproof

We turn to the second proof of Theorem~\ref{thm_stir}, by so called Laplace's method, and we follow N.\,G. de Bruijn \cite[Chapter 4]{debr}. We start with a classical formula, due essentially but not entirely to L. Euler. By V.\,S. Varadarajan \cite[p. 100]{vara}, 
L. Euler would write the gamma function integral below as 
$$
\int_0^1(-\log x)^n\;dx\;,
$$
and it was A.-M. Legendre who wrote it in the familiar form in the infinite range. Our next calculation is less anachronistic than some of the others 
because in the times of L. Euler and A.-M. Legendre there were only Newton integrals. The second proof uses integrals over infinite intervals and for their existence we cannot rely on Corollary~\ref{prop_ex_newt}.

\begin{prop}[second expression of $n!$ by an $\int$]\tec\label{gamma}
For all $n$ in $\N_0$ we have
$$
n!=(N)\int_0^{+\infty}x^n e^{-x}=(N)\int_0^{+\infty}x^n e^{-x}\,dx\;.
$$
\end{prop}
\proof
We denote the integral by $I_n$. We prove its existence and compute its value by induction on $n$. First, $I_0=(-e^{-x})(+\infty^-)-(-e^{-x})(0^+)=0-(-1)=1$. For $n>0$ we get by Theorem~\ref{prop_pp_newt}
(in the integration by parts identity below the second term is clearly defined and finite, and so is the third by the inductive assumption) and by part 1 of Proposition~\ref{prop_linear} that
\begin{eqnarray*}
I_n&=&(N)\int_0^{+\infty}x^n(-e^{-x})'=(-x^ne^{-x})(+\infty^-)-(-x^ne^{-x})(0^+)+\\
&&+\;(N)\int_0^{+\infty}(x^n)'e^{-x}=0-0+n\cdot(N)\int_0^{+\infty}x^{n-1}e^{-x}\\
&=&nI_{n-1}\;.
\end{eqnarray*}
By induction, $I_n$ exists for every $n\in\N_0$ and $I_n=n!$.
\eproof

Let $n\in\N$. Substitution $x\leftarrow y$, $x=n(1+y)$, by Proposition~\ref{prop_subst} gives
$$
(N)\int_0^{+\infty}x^n e^{-x}=e^{-n}n^{n+1}\cdot(N)\int_{-1}^{+\infty}(e^{-y}(1+y))^n\;.
$$
Let $f(y)=e^{-y}(1+y)$. Then $f'(y)=-e^{-y}y>0$ on $[-1,0)$ and is $<0$ on $(0,+\infty)$, and we see that $f(y)$ increases from $0$ to $1$ on $[-1,0]$ and decreases from 
$1$ to $0^+$ on $[0,+\infty)$. We identify intervals around $0$ with the bulk of the last integral concentrated in them, and replace the integrand with a neater function.

\begin{prop}[concentration of the $\int$]\tec\label{prop_bulk_conce}
If $\de=\de(n)\cc\N\to(0,1)$ is a~sequence such that $n\de^3\to0$
as $n\to\infty$, then for all $n\in\N$ one has
$$
(N)\int_{-1}^{+\infty}(e^{-y}(1+y))^n=(1+O(n\de^3))\cdot(N)\int_{-\de}^{\de}e^{-ny^2/2}+O(e^{-n\de^2/2})\;.
$$
\end{prop}
\proof
Using again the expansion of $\log(1+x)$ we have 
\begin{eqnarray*}
f(y)=e^{-y}(1+y)&=&\exp(-y+\log(1+y))=\exp(-y^2/2+O(y^3))\\
&=&e^{-y^2/2}(1+O(y^3))\ (y\in[-1/2,\,2],\mbox{ say})\;.
\end{eqnarray*}
If $\de=\de(n)$ is as stated then 
$$
(1+O(\de^3))^n=\exp(n\log(1+O(\de^3)))=\exp(O(n\de^3))=1+O(n\de^3)
$$
and $f(-\de)^n, f(\de)^n=O(e^{-n\de^2/2})$. Using part 2 of Proposition~\ref{prop_linear} we define the decomposition
\begin{eqnarray*}
(N)\int_{-1}^{+\infty}f(y)^n&=&(N)\int_{-1}^{-\de}+\;(N)\int_{-\de}^{\de}+\;(N)\int_{\de}^4+\;(N)\int_4^{+\infty}\\
&=:&I_1+I_2+I_3+I_4\;.
\end{eqnarray*}
Since $0\le f(y)^n\le f(-\de)^n$ on $[-1,-\de]$ and $0< f(y)^n\le f(\de)^n$ on $[\de,+\infty)$, the above estimates and Corollary~\ref{ML_bound} show that both 
$I_1,I_3=O(e^{-n\de^2/2})$. Since $1+y\le 1+y/2+y^2/8\le e^{y/2}$ for $y\ge4$, $f(y)\le e^{-y/2}$ for $y\ge4$ and by Proposition~\ref{prop_mono_newt},
$$
I_4\le(N)\int_4^{+\infty}e^{-ny/2}=\frac{2e^{-2n}}{n}\;.
$$
So $I_4=O(e^{-n\de^2/2})$ too. The remaining integral satisfies
\begin{eqnarray*}
I_2&=&(N)\int_{-\de}^{+\de}f(y)^n=(N)\int_{-\de}^{+\de}e^{-ny^2/2}(1+O(ny^3))\\
&=&(1+O(n\de^3))\cdot(N)\int_{-\de}^{\de}e^{-ny^2/2}
\end{eqnarray*}
(the last equality follows by Proposition~\ref{prop_mono_newt}) and we are done. 
\eproof

\begin{prop}[reduction to the Gauss $\int$]\tec\label{redu_to_gauss}
If $\de=\de(n)$ is as in the previous proposition and $m\in\N$ then
$$
(N)\int_{-\de}^{\de}e^{-ny^2/2}=\sqrt{\frac{2}{n}}\cdot(N)\int_{-\infty}^{+\infty}e^{-t^2}+O(e^{-n\de^2/2})\;.
$$
\end{prop}
\proof
We define, using part 2 of Proposition~\ref{prop_linear}, eveness of the integrand, and the version of Proposition~\ref{prop_subst} with the
flipping substitution $g(y)=-y$, the decomposition
\begin{eqnarray*}
(N)\int_{-\de}^{\de}e^{-ny^2/2}&=&(N)\int_{-\infty}^{+\infty}e^{-ny^2/2}-2\cdot(N)\int_{\de}^{+\infty}e^{-ny^2/2}\\
&=:&I_5-2I_6\;,
\end{eqnarray*}
provided that $I_5$ exists. We are in a similar situation as in the beginning of the proof of Corollary~\ref{suma}. But $I_5$ exists by the majorization $e^{-a^2}\le e^{-a}$ for $a\ge1$ (as we already know from the beginning of the proof of Theorem~\ref{fubi_infi_inter}). We estimate $I_6$ in the same way as we estimated $I_3$ and $I_4$ in the previous proof and get 
the same bound $I_6=O(e^{-n\de^2/2})$. Proposition~\ref{prop_subst} with the substitution 
$y\leftarrow t$, $y=t\sqrt{2/n}$, yields
$$
I_5=\sqrt{\frac{2}{n}}\cdot(N)\int_{-\infty}^{+\infty}e^{-t^2}\;.
$$
\eproof

It remains to compute the Gauss integral $(N)\int_{-\infty}^{+\infty}e^{-t^2}$ and to select the sequence $\de=\de(n)$.

\begin{prop}[the Gauss $\int$]\tec\label{gauss_int}
We have the identity
$$
(N)\int_{-\infty}^{+\infty}e^{-t^2}=\sqrt{\pi}\;.
$$
\end{prop}
\proof
By part 2 of Proposition~\ref{prop_linear} and the version of Proposition~\ref{prop_subst} with the flipping substitution $g(t)=-t$, we need to prove that
$$
I_7:=(N)\int_0^{+\infty}e^{-t^2}=\frac{\sqrt{\pi}}{2}
$$
(in the beginning of the proof of Theorem~\ref{fubi_infi_inter} we proved that $I_7$ exists). This is equivalent with $I_7^2=\pi/4$. Indeed, we compute (we justify each of the eight 
steps after the computation)
\begin{eqnarray*}
I_7^2&=&(N)\int_0^{+\infty}e^{-t^2}\cdot(N)\int_0^{+\infty}e^{-u^2}=(N)\int_0^{+\infty}\left((N)\int_0^{+\infty}e^{-t^2}\right)e^{-u^2}\\
&=&(N)\int_0^{+\infty}\left((N)\int_0^{+\infty}e^{-t^2-u^2}\;dt\right)\;du\\
&=&(N)\int_0^{+\infty}\left((N)\int_0^{+\infty}ue^{-u^2(1+v^2)}\;dv\right)\;du\\
&=&(N)\int_0^{+\infty}\left((N)\int_0^{+\infty}ue^{-u^2(1+v^2)}\;du\right)\;dv=(N)\int_0^{+\infty}\frac{1}{2(1+v^2)}\\
&=&\frac{\arctan(+\infty^-)-\arctan(0^+)}{2}=\frac{\pi}{4}\;.
\end{eqnarray*}
The first four steps repeat the computation from the beginning of the proof of Theorem~\ref{fubi_infi_inter}, and the crucial fifth step is this theorem.
In the sixth step we compute the inner integral according to Definition~\ref{defi_newt} by the  primitive (for fixed $v\in\R$)
$$
\frac{d}{du}\left(\frac{-e^{-u^2(1+v^2)}}{2(1+v^2)}\right)=ue^{-u^2(1+v^2)}\;.
$$
In the last two steps we compute the integral according to Definition~\ref{defi_newt} by the primitive $(\arctan v)'=\frac{1}{1+v^2}$.
\eproof

\noindent
The number $\pi/2$ came about again as the smallest positive root of $\cos x$. Computation of the Gauss integral just by
the Newton integration may be of some interest, for in B. Conrad \cite[p. 1]{conr} we read: 
``$\Phi(u)=\frac{1}{\sqrt{2\pi}}\int_{-\infty}^ue^{-u^2/2}\mathrm{d}u$ (\ds) so the evaluation of $\Phi(\infty)$ must proceed by a
method different from the calculation of anti-derivatives as in calculus.'' The point of our article is that anti-derivatives 
(primitives) fully suffice for such evaluation. But the difficulties with Theorem~\ref{fubi_infi_inter} show that it is not as straightforward as one might think.

To finish, we set $\de=\de(n)=n^{-1/2+\ep/3}$ where $\ep\in(0,1/2)$. Combining Propositions~\ref{gamma}--\ref{gauss_int}
we obtain the asymptotics
\begin{eqnarray*}
n!&=&e^{-n}n^{n+1}\left((1+O(n\de^3)).(\sqrt{2\pi/n}+O(e^{-n\de^2/2}))+O(e^{-n\de^2/2})\right)\\
&=&\sqrt{2\pi n}\left(\frac{n}{e}\right)^n(1+O(n^{-1/2+\ep}))
\end{eqnarray*}
because $O(e^{-n\de^2/2})=O(e^{-n^{2\ep/3}/2})$ goes to $0$ for $n\to\infty$ faster than $n^{-c}$ for any $c>0$. This completes the second proof of Theorem~\ref{thm_stir}. \eproof

\bigskip\noindent
{\em Martin Klazar\\
Department of Applied Mathematics\\
Charles University, Faculty of Mathematics and Physics\\
Malostransk\'e n\'am\v est\'\i\ 25\\
11800 Praha\\
Czechia\\
{\tt klazar@kam.mff.cuni.cz}
}
\end{document}